\documentclass[11pt,a4]{article}
\usepackage[T2A]{fontenc}
\usepackage[cp866]{inputenc}
\usepackage[mathscr]{eucal}
\usepackage{amssymb}
\usepackage{amsmath}
\renewcommand{\leq}{\leqslant}

\renewcommand{\C}{{\mathbb C}}

\newcommand{\R}{{\mathbb R}}

\newcommand{\K}{{\mathbb K}}

\usepackage{fancyhdr}
\begin{document}
\sloppy

{\normalsize

\thispagestyle{empty}

\mbox{}
\\[-2.75ex]
\centerline{
\large
\bf
INVERSE PROBLEM OF GROUP ANALYSIS FOR 
}
\\[0.75ex]
\mbox{}\hfill
{
\large
\bf
AUTONOMOUS DIFFERENTIAL SYSTEMS\;\!}\footnote[1]{
The basic results of this paper have been published in the monograph
{\it {\rm"}Integrals of differential systems}", Grodno, 2006, 447 p. [1, pp. 152 -- 160] and
in the journals [2 -- 5]: 
{\it Differential Equations}, Vol. 30, 1994, No. 5, 899-901; 
{\it Vestsi Nats. Akad. Navuk Belarusi}, Ser. 2, 1995, No. 2, 46-49;   
{\it Proceedings of the scientific seminars POMI}, Vol. 232, 1996, 86-89;
{\it Dokl. Akad. Nauk Belarusi}, Vol. 50, 2006, No. 3,  15-19.
}
\hfill\mbox{}
\\[2.5ex]
\centerline{
\bf 
V.N. Gorbuzov
}
\\[1.75ex]
\centerline{ 
\it 
Department of mathematical analysis and differential equations
} 
\\[1.25ex]
\centerline{ 
\it 
Faculty of Mathematics and Computer Science
} 
\\[1.25ex]
\centerline{ 
\it 
Yanka Kupala Grodno State University
} 
\\[1.25ex]
\centerline{ 
\it
Ozeshko 22, Grodno, 230023, Belarus
}
\\[1.25ex]
\centerline{ 
E-mail: gorbuzov@grsu.by
}
\\[5ex]
\centerline{\large\bf Abstract}
\\[1ex]
\indent
Full set of autonomous completely solvable differential systems of equations in total differentials is 
built by basis of infinitesimal operators, universal invariant, and structure constants of  admited
multiparametric Lie group (abelian and non-abelian).
\\[1.5ex]
\\[1.25ex]
%58A17 Pfaffian systems
%35F35 Linear first-order systems
%34A34 Nonlinear equations and systems (ordinary), general
%35R03 Partial differential equations on Heisenberg groups, Lie groups, Carnot groups, etc.
%57S15 Compact Lie groups of differentiable transformations
%57S20 Noncompact Lie groups of transformations
%70G65 Symmetries, Lie-group and Lie-algebra methods (mechanics)
%22E05 Local Lie groups
%22E30 Analysis on real and complex Lie groups
\\[5ex]
\indent
{\bf Problem statement.}
Consider an autonomous completely solvable system of equations in total differentials
\\[1.75ex]
\mbox{}\hfill                                   %(1)
$
dx=f(x)\, dt,  \quad
dx={\rm colon}\,(dx_1^{},\ldots,dx_n^{}), \ \
dt={\rm colon}\,(dt_1^{},\ldots,dt_m^{}),
$
\hfill (1)
\\[2ex]
on a domain $\Omega$ from the phase space $\K^n$  (real $\R^n$ or complex $\C^n)$ 
\vspace{0.35ex}
and let the ent\-ries of the $(n\times m)\!$-matrix  
\vspace{0.35ex}
$f(x)=\bigl\| f_{ij}^{}(x)\bigr\|$ for all $x\in \Omega$ 
be holomorphic functions $f_{ij}^{}\colon \Omega\to \K,\ i=1,\ldots,n,\ j=1,\ldots,m.$
\vspace{0.5ex}

The linear differential ope\-ra\-tors of first order
\\[1.75ex]
\mbox{}\hfill                                  
$
\displaystyle
{\frak F}_j^{}(x)=\sum_{i=1}^n\, f_{ij}^{}(x)\;\!\partial_{x_i^{}}^{}
$ 
\ for all 
$
x\in \Omega, 
\quad
j=1,\ldots,m,
\hfill
$
\\[2ex]
are autonomous operators of differentiation by virtue of sys\-tem (1) [1, p. 20; 6; 7].

An operational criterion of completely solvability on the domain $\Omega$ of system (1) is represented via 
Poisson brackets as the system of  identities [1, pp. 17 -- 25; 6]
\\[2ex]
\mbox{}\hfill                                  
$
\displaystyle
\bigl[ {\frak F}_j^{}(x), {\frak F}_\zeta^{}(x) \bigr]={\frak O}
$
\ for all 
$
x\in \Omega, 
\quad 
j=1,\ldots,m,\ \zeta=1,\ldots,m.
\hfill
$
\\[2ex]
 \indent
Suppose 
\vspace{0.25ex}
a $q\!$-parameter Lie group of transformation $G_q^{},\, 1\leq q\leq n,$ [8 -- 11]
has $(n-k)\!$-cy\-lindrical infinitesimal operators 
\vspace{0.25ex}
(coordinates of these operators not depends on $n-k$ dependent variables [1, pp. 106 -- 167;  7])
\\[2ex]
\mbox{}\hfill               % (2)
$
\displaystyle
{\frak G}_{l}^{}(x) =
\sum_{i=1}^{n}\,
{\rm g}_{li}^{}({}^k x)\;\!\partial_{x_i^{}}^{}
$
\ for all 
$
x \in \Omega, 
\ \ l = 1,\ldots,q,
\quad
{}^k x=(x_1^{},\ldots,x_k^{}),\ 1\leq k\leq n,
$
\hfill (2)
\\[2ex]
and a basis of absolute invariants
\\[1.75ex]
\mbox{}\hfill
$
I\colon x\to  \bigl(I_1^{}(x), \ldots, I_{n-q}^{}(x)\bigr)
$
\ for all 
$
x \in \Omega.
$
\hfill (3)
\\[2ex]
\indent
Besides, the Poisson brackets 
\\[2ex]
\mbox{}\hfill                           % (4)
$
\displaystyle
\bigl[{\frak G}_l^{}(x), {\frak G}_s^{}(x)\bigr] =
\sum_{p=1}^{q}\,c_{lsp}^{}\;\!
{\frak G}_p^{}(x)
$ 
\ for all 
$
x \in  \Omega, 
\quad
l= 1,\ldots, q,
\ \ s = 1,\ldots, q,
$
\hfill (4)
\\[2.25ex]
where
\vspace{0.75ex}
$\!c_{lsp}^{},\, l\!=\! 1,\ldots, q,\, s\! =\! 1,\ldots, q,\, p\! =\! 1,\ldots, q,\!$ 
are the structure constants of Lie group $\!G_q^{}.$

Assume that 
\vspace{0.25ex}
the coordinate functions ${\rm g}_{li}^{},\ l = 1,\ldots, q,\ i =1,\ldots, n,$  
of the infinitesimal operators (2) and the absolute invariants $I_{\tau}^{},\ \tau= 1,\ldots, n-q,$ of the basis (3) 
are holomorphic on the domain $\Omega.$ Furthermore, the infinitesimal operators (2)
are not holomorphically linearly bound [1, p. 11; 6] on the domain $\Omega.$

In this paper, the following inverse problem of group analysis for differential systems is considered:
to select from the set off all completely solvable differential systems (1) the systems  that admit a Lie group $G_q^{}$ 
with the infinitesimal operators (2), the universal invariant (3), and the structure constants $c_{lsp}^{}$ from 
the operator identities (4).
\vspace{0.5ex}

{\bf  
Classes of differential systems that admit an non-abelian Lie group.}
Suppose the completely solvable autonomous system of equations in total differentials (1) admits 
a Lie group $G_q^{}.$ Then, the Poisson brackets
\\[2ex]
\mbox{}\hfill                           % (5)
$
\displaystyle
\bigl[{\frak G}_l^{}(x), {\frak F}_j^{}(x)\bigr]=
{\frak O} 
$
\ for all 
$
x \in  \Omega, 
\quad
l = 1,\ldots, q, \ \ j= 1,\ldots, m.
$
\hfill (5)
\\[2.5ex]
\indent
{\bf Lemma.}
{\it
Suppose the linear differential ope\-ra\-tors of first order {\rm (2)} are 
not holomorphically linearly bound on the domain $\Omega.$
Then the following conditions are true{\rm:}

{\rm i)} 
the complete linear homogeneous system of partial differential equations
}
\\[2ex]
\mbox{}\hfill                   %  (6)
$
{\frak G}_{l}^{}(x)\;\! u = 0, 
\quad 
l = 1,\ldots, q,  
\quad 
q < n,
$
\hfill (6)
\\[2ex]
{\it
admits $(n-q)\!$-parameter abelian Lie group $\widetilde{G}_{n-q}^{};$

{\rm ii)} infinitesimal operators of the abelian Lie group $\widetilde{G}_{n-q}^{}$ and
\vspace{0.35ex}
the operators $\!{\frak G}_l^{},\, l\! =\! 1,\ldots, q,$ are commutative and they are 
not holomorphically linearly bound on the domain $\Omega.$
}
\vspace{0.5ex}

{\sl Proof}
(analogous statement for linear homogeneous partial differential equation was proved in [9, p. 109]).
\vspace{0.25ex}
Since the system (6) is complete [1, p. 38; 6; 12, p. 70], we see that 
integral basis of this system is $n - q$ first integrals [1, p. 50]  
\\[2ex]
\mbox{}\hfill
$
F_\xi^{}\colon x\to F_\xi^{}(x)
$
\ for all 
$
x\in \Omega,
\quad
\xi = 1,\ldots, n-q.
\hfill
$
\\[2ex]
\indent
Then, the substitution
\\[2ex]
\mbox{}\hfill
$
\nu_\xi^{} = F_\xi^{}(x), 
\ \ \xi = 1,\ldots, n-q, 
\quad \
\nu_\zeta^{} =\Phi_\zeta^{}(x), 
\ \ \zeta =n-q+1,\ldots, n,
\hfill
$
\\[2.25ex]
where the functions 
\vspace{0.5ex}
$\Phi_\zeta^{},\ \zeta =n-q+1,\ldots, n,$ 
are first integrals of the 
complete linear non-homogeneous partial differential systems [12, p. 101]
\\[2ex]
\mbox{}\hfill
$
{\frak G}_l^{}(x)\;\! u =
\delta_{l, \zeta-n+q}^{}\;\!, 
\quad
l =  1,\ldots, q, \ \ 
\zeta = n - q +1,\ldots, n, \ \ 
(\delta_{ls}^{}\ \text{is Kronecker symbol}),
\hfill
$
\\[2.25ex]
reduces the system (6) to the 
complete linear homogeneous partial differential system
\\[2ex]
\mbox{}\hfill
$
\partial_{\nu_\zeta^{} }^{}\;\! u = 0, 
\quad 
\zeta = n - q +1,\ldots, n.
\hfill
$
\\[2ex]
\indent
This system admits the $n - q$ differential 
operators $\partial_{\nu_\xi^{}}^{}, \ \xi = 1,\ldots, n-q.$ 
At the same time these operators
are not holomorphically linearly bound on the space $\K^{n-q}$ and they admits
$(n-q)\!$-parameter abelian Lie group.
\vspace{0.25ex}

Consequently the system (6) admits 
$(n-q)\!$-parameter abelian Lie group $\widetilde{G}_{n-q}^{}.$
\vspace{0.35ex}

The operators $\partial_{\nu_i^{}}^{}, \ i = 1,\ldots, n,$ are commutative and these operators
\vspace{0.25ex}
are not holomorphically linearly bound on the space $\K^{n}.$ 
In addition, 
the property of symmetry for Poisson brackets and 
the property of linear bound (or not linear bound) for operators
are conservation under transformation. This implies that:
\vspace{0.35ex}
infinitesimal operators of the abelian Lie group $\widetilde{G}_{n-q}^{}$ and
\vspace{0.35ex}
the operators ${\frak G}_l^{},\ l = 1,\ldots, q,$ are commutative;
infinitesimal operators of the abelian Lie group $\widetilde{G}_{n-q}^{}$ and
\vspace{0.5ex}
the operators ${\frak G}_l^{},\ l = 1,\ldots, q,$ are not holomorphically linearly bound on the domain $\Omega.$
The lemma is proved.
\vspace{0.75ex}

Using the operators (2) admitted by the system (1), we can build the
$(n-k)\!$-cy\-lindrical [7] linear homogeneous system of partial differential equations 
\\[2ex]
\mbox{}\hfill                                      % (7)
$
{\frak G}_l^{}(x)\;\!u =  0, 
\quad
l = 1,\ldots,q, 
\quad
q \leq n.
$
\hfill (7)
\\[2ex]
\indent
By the conditions (4), it follows that the system (7) is complete [1, p. 38; 6] on $\Omega.$
\vspace{0.35ex}

From Lemma it follows that the complete system (7)  admits an
\vspace{0.25ex}
$(n-q)\!$-parameter abelian group $\widetilde{G}_{n-q}^{}.$
Let the linear differential ope\-ra\-tors of first order
\\[2ex]
\mbox{}\hfill
$
\displaystyle
{\frak G}_\tau(x) =
\sum\limits_{i=1}^{n}\,{\rm g}_{\tau i}^{}(x)\;\!\partial_{x_i^{}}^{}
$
\ for all 
$
x \in \Omega, 
\quad 
\tau = q+1,\ldots, n,
\hfill
$
\\[2ex]
be infinitesimal operators of the abelian Lie group $\widetilde{G}_{n-q}^{}.$
\vspace{0.35ex}

By Lemma, the Poisson brackets
\\[2.25ex]
\mbox{}\hfill                                    %(8)
$
\bigl[{\frak G}_i^{}(x), {\frak G}_\tau^{}(x)\bigr]  =  {\frak O}
$
\ for all 
$
x \in \Omega,
\quad
i =1,\ldots, n, \ \ \tau = q+1,\ldots, n.
$
\hfill (8)
\\[2.26ex]
Note also that the linear differential operators ${\frak G}_i^{}, \ i =1,\ldots, n,$ on the domain  $\Omega$
are not holomorphically linearly bound. 

Let the operators ${\frak G}_i^{}, \ i =1,\ldots, n,$ be a basis.
\vspace{0.25ex}
Then the operators ${\frak F}_j^{}, \ j = 1,\ldots, m,$ induced by the system (1) have the expansions
\\[2ex]
\mbox{}\hfill                                         % (9)
$
\displaystyle
{\frak F}_{j}^{}(x)
=
\sum_{i=1}^{n}\,\psi_{ji}^{}(x)\;\! {\frak G}_{i}^{}(x)
$
\ for all 
$
x \in \Omega, \ \  j = 1,\ldots, m,
$
\hfill (9)
\\[2.35ex]
where  $\psi_{ji}^{}\colon \Omega\to \K,\ j = 1,\ldots, m,\ i =1,\ldots, n,$ are holomorphic functions.
\vspace{0.75ex}

Using the expansions (9) and the commutator identities (5) and (8), we obtain 
necessary and sufficient conditions that the completely solvable system (1) admits the Lie group $G_q^{}\colon$
\\[2ex]
\mbox{}\hfill
$
\displaystyle
\bigl[{\frak G}_{l}^{}(x), {\frak F}_{j}^{}(x)\bigr] = 
\sum_{i=1}^{n}\,
{\frak G}_{l}^{}\psi_{ji}^{}(x)\;\!
{\frak G}_{i}^{}(x)  +
\sum_{i=1}^{q}\;\!\sum_{p=1}^{q}\,
c_{lpi}^{}\;\! \psi_{jp}^{}(x)\;\!
{\frak G}_{i}^{}(x)
$
\ for all 
$
x \in \Omega, 
\hfill
$
\\[2ex]
\mbox{}\hfill
$
l = 1,\ldots, q, 
\quad
j =1,\ldots, m.
\hfill
$
\\[2.5ex]
\indent
Since the commutator identities (5) hold and 
the operators ${\frak G}_i^{}, \ i = 1, \ldots,n,$ are basis, we have
\\[2ex]
\mbox{}\hfill                                 % (10)
$
{\frak G}_{l}^{}\;\! \psi_{j\theta}^{}(x)
\displaystyle
+ \sum_{p=1}^{q}c_{lp\theta}^{}\;\!
\psi_{j p}^{}(x) = 0
$ 
\ for all 
$
x \in \Omega, 
\quad  
l = 1,\ldots, q, \ \
j =1,\ldots, m, \ \
\theta = 1,\ldots, q, 
\hfill
$
\\[0ex]
\mbox{}\hfill  (10)
\\
\mbox{}\hfill
$
{\frak G}_{l}^{}\;\! \psi_{j\tau}^{}(x) = 0
$
\ for all 
$
x \in \Omega, 
\quad
l = 1,\ldots, q, \ \
j =1,\ldots, m, \ \
\tau  = q+1,\ldots, n.
\hfill
$
\\[2.5ex]
\indent
From this system of identities it follows that the scalar functions
\vspace{0.25ex}
$
\psi_{j\tau}^{}, \
j =1,\ldots, m, 
\linebreak
\tau  = q+1,\ldots, n,
$
are absolute invariants of Lie group $G_q^{}.$ Therefore,
\\[2ex]
\mbox{}\hfill
$
\psi_{j\tau}^{}(x)=
\varphi_{j\tau}^{}(I(x))
$
\ for all 
$
x \in \Omega,  
\quad
j =1,\ldots, m, \ \
\tau  = q+1,\ldots, n,
\hfill
$
\\[2.35ex]
where the functions
\vspace{0.55ex}
$\varphi_{j\tau}^{}\colon  W \to \K, \ 
j =1,\ldots, m, \
\tau  = q+1,\ldots, n,$ 
are holomorphic on the domain $W$ of the space $\K^{n-q}.$

Taking into account the system of identities (10), 
\vspace{0.35ex}
we obtain the scalar functions
$
\psi_{j\theta}^{},
\linebreak
j =1,\ldots, m, \
\theta = 1,\ldots, q, 
$
\vspace{0.25ex}
are first integrals of the 
linear non-homogeneous system of partial differential equations 
\\[2ex]
\mbox{}\hfill                                 % (11)
$
{\frak G}_{l}^{}(x)\;\!\psi_{j\theta}^{}
\displaystyle
+ \sum_{p=1}^{q}c_{lp\theta}^{}\;\!
\psi_{j p}^{}= 0, 
\quad
l = 1,\ldots, q, \ \
j =1,\ldots, m, \ \
\theta = 1,\ldots, q.
$
\hfill (11)
\\[2.25ex]
\indent
By the identities (4) and  (8), it follows that the system (11) is complete on the domain $\Omega.$ 
Suppose this system has the first integrals
\\[2ex]
\mbox{}\hfill
$
\psi_{jl}^{}\colon x\to\
\psi_{jl}^{}(x)
$
\ for all 
$
x\in \Omega, 
\quad
j =1,\ldots, m, \ \
l = 1,\ldots, q.
\hfill
$
\\[2.2ex]
\indent
Then, we obtain the representation
\\[2ex]
\mbox{}\hfill
$
\displaystyle
{\frak F}_{j}^{}(x)  =
\sum_{l=1}^{q}\psi_{j l}^{}(x)\;\! {\frak G}_{l}^{}(x)  +
\sum_{\tau=q+1}^{n}\varphi_{j\tau}^{}(I(x))\;\!
{\frak G}_{\tau}^{}(x)
$
\ for all  
$
x\in \Omega, 
\quad
j =1,\ldots, m, 
\hfill
$
\\[2.25ex]
and the following statement.
\vspace{1ex}

{\bf Theorem 1.}
\vspace{0.15ex}
{\it
The autonomous completely solvable system of equations in total differentials {\rm (1)} admits 
the $q\!\!$-parameter Lie group $G_q^{}$ with
\vspace{0.15ex}
$(n-k)\!$-cy\-lindrical infinitesimal operators {\rm (2),}
the universal invariant {\rm (3)}, and the structure constants from the representations {\rm(4)}
\vspace{0.15ex}
if and only if this differential system has the form}
\\[2.25ex]
\mbox{}\hfill          % (12)
$
dx_i^{} = \displaystyle
\sum_{j=1}^{m}\;\!
\biggl(\,
\sum_{l=1}^{q}\psi_{jl}^{}(x)\;\!
{\rm g}_{li}^{}({}^k\!x) +
\sum_{\tau = q+1}^{n}
\varphi_{j \tau}^{}(I(x))\;\!
{\rm g}_{\tau i}^{}(x)\biggr)\;\! dt_j^{},
\quad 
i = 1,\ldots, n,
$
\hfill (12)
\\[2.5ex]
{\it
where the holomorphic functions
\vspace{0.35ex}
$\psi_{jl}^{}\colon
\Omega \to \K,\ j = 1,\ldots, m,\
l = 1, \ldots,q,$ are first integrals of the
\vspace{0.15ex}
complete linear non-homogeneous system of partial differential equations {\rm (11),} 
the holomorphic functions 
\vspace{0.35ex}
${\rm g}_{\tau i}^{}\colon  \Omega\to \K,\ \tau= q+1,\ldots, n,\ i =1,\ldots, n,$ are 
the coordinates of the infinitesimal operators 
\vspace{0.5ex}
${\frak G}_\tau^{},\, \tau =q+1,\ldots, n,$ of
$(n-q)\!\!$-parameter abelian Lie group $\widetilde{G}_{n-q}^{}$ admitted by the  
\vspace{0.35ex}
$(n-k)\!$-cy\-lindrical complete linear homogeneous system of partial differential equations {\rm (7),} 
\vspace{0.35ex}
and the holomorphic on the domain $W\subset \K^{n-q}$ functions 
$\varphi_{j\tau}^{}\colon W \to \K,\ j = 1,\ldots, m,\ \tau = q+1,\ldots,n,$ such that
\vspace{0.5ex}
the conditions of completely solvability on the domain $\Omega$ for system {\rm (12)} hold.
} 
\vspace{1.5ex}

{\bf  
Classes of differential systems that admit an abelian Lie group.}
\vspace{0.15ex}
Suppose the $q\!$-parameter Lie group $G_q^{}$ is abelian. 
Then in the representation (4) the structure constants 
\\[2ex]
\mbox{}\hfill
$
c_{lsp}^{}=0,
\quad 
l=1,\ldots,q,\ \ 
s=1,\ldots,q,\ \ 
p=1,\ldots,q.
\hfill
$
\\[2.25ex]
\indent
Also, from the system of identities (10) it follows that the functions 
\vspace{0.5ex}
$\psi_{ji}^{},\ j=1,\ldots,m,$ $i=1,\ldots,n,$ 
are absolute invariants of the $q\!$-parameter abelian Lie group $G_q^{}.$ Therefore,
\\[2ex]
\mbox{}\hfill
$
\psi_{ji}^{}(x)=\varphi_{ji}^{}(I(x))
$
\ for all 
$
x\in\Omega,
\quad 
j=1,\ldots,m,\ \ i=1,\ldots,n,
\hfill
$
\\[2.25ex]
where the functions 
\vspace{0.5ex}
$\varphi_{ji}^{}\colon W\to\K,\ j=1,\ldots,m,\ i=1,\ldots,n,$ 
are holomorphic on the domain $W$ from the space $\K^{n-q}.$
\vspace{0.35ex}

Then, 
\\[2ex]
\mbox{}\hfill
$
\displaystyle
{\frak F}_{j}^{}(x)  =
\sum_{l=1}^{q}\,\varphi_{j l}^{}(I(x))\;\!{\frak G}_{l}^{}(x)  +
\sum_{\tau=q+1}^{n}\,\varphi_{j\tau}^{}(I(x))\;\!
{\frak G}_{\tau}^{}(x)
$
\ for all  
$
x\in \Omega, 
\quad
j =1,\ldots, m, 
\hfill
$
\\[2.2ex]
and we have

{\bf Theorem 2.}
{\it
The autonomous completely solvable system of equations in total differentials {\rm (1)} admits 
the $q\!\!$-parameter abelian Lie group $G_q^{}$ with
$(n-k)\!\!$-cy\-lindrical infinitesimal operators {\rm (2)} and
the universal invariant {\rm (3)} if and only if this system has the form}
\\[2.25ex]
\mbox{}\hfill        
$
\displaystyle
dx_i^{} = 
\sum_{j=1}^{m}\;\!
\biggl(\,
\sum_{l=1}^{q} \varphi_{jl}^{}(I(x))\;\!
{\rm g}_{li}^{}({}^k\!x) +
\sum_{\tau = q+1}^{n}
\varphi_{j \tau}^{}(I(x))\;\!
{\rm g}_{\tau i}^{}(x)\biggr)\;\! dt_j^{},
\quad 
i = 1,\ldots, n,
\hfill
$
\\[2.5ex]
{\it
where the holomorphic functions 
\vspace{0.25ex}
${\rm g}_{\tau i}^{}\colon  \Omega\to \K,\ \tau= q+1,\ldots, n,\ i =1,\ldots, n,$ are 
the coordinates of the infinitesimal operators 
\vspace{0.25ex}
${\frak G}_\tau^{},\, \tau =q+1,\ldots, n,$ of
$(n-q)\!\!$-parameter abelian Lie group $\widetilde{G}_{n-q}^{}$ admitted by the  
\vspace{0.25ex}
$(n-k)\!$-cy\-lindrical complete linear homogeneous system of partial differential equations {\rm (7),} 
\vspace{0.25ex}
and the holomorphic on the domain $W\subset \K^{n-q}$ functions 
$\varphi_{j\;\!i}^{}\colon W \to \K,\ j = 1,\ldots, m,\ i =1,\ldots,n,$ such that
\vspace{0.25ex}
the conditions of completely solvability on the domain $\Omega$ for this built differential system hold.
}
\vspace{1ex}

{\bf Example 1.}
Consider the $q\!$-parameter abelian Lie group of transformations
\\[2ex]
\mbox{}\hfill          %(13)
$
x_l^{}\to\, e^{\;\!\alpha_l^{}}\;\!x_l^{},
\ \ l=1,\ldots,q,
\quad 
x_{\tau}^{}\to x_{\tau}^{},
\ \ \tau=q+1,\ldots,n,
\ \ q<n,
$
\hfill (13)
\\[2.25ex]
with $(n-q)\!$-cy\-lindrical infinitesimal operators
\\[2ex]
\mbox{}\hfill            %(14)
$
{\frak G}_l^{}(x)=x_l^{}\;\!\partial_{x_l^{}}^{}
$
\ for all 
$
x\in\K^n,
\quad 
l=1,\ldots,q,
$
\hfill (14)
\\[2.25ex]
and the basis of absolute invariants
\\[2ex]
\mbox{}\hfill             %(15)
$
I\colon x\to (x_{q+1}^{},\ldots, x_n^{})
$
\ for all 
$
x\in\K^n.
$
\hfill (15)
\\[2.25ex]
\indent
By Theorem 2, we have the following
\vspace{0.75ex}

{\bf Proposition 1.}
{\it 
The completely solvable autonomous system of equations in total differentials {\rm (1)} admits 
the $q\!\!$-parameter abelian Lie group of transformations {\rm (13)} with the
infinitesimal operators {\rm (14)} and
the universal invariant {\rm (15)} if and only if this system has the form
}
\\[2.25ex]
\mbox{}\hfill        %(16)
$
\displaystyle
dx_i^{} = 
\sum_{j=1}^{m}\;\!
\biggl(\,
\sum_{l=1}^{q}\! \delta_{il}^{}\;\! x_l^{}\;\! \varphi_{jl}^{}(x_{q+1}^{},\ldots,x_{n}^{})
+
\sum_{\tau = q+1}^{n}
{\rm g}_{\tau i}^{}(x)\;\! \varphi_{j \tau}^{}(x_{q+1}^{},\ldots,x_{n}^{}  )\!\biggr) dt_j^{},
\ i\! =\! 1,\ldots, n,
$
\hfill (16)
\\[2.5ex]
{\it
where the holomorphic functions 
${\rm g}_{\tau i}^{}\colon  \Omega\to \K,\ 
\tau= q+1,\ldots, n,\ i =1,\ldots, n,$
\vspace{0.5ex}
are the coor\-di\-na\-tes of the linear differential ope\-ra\-tors of first order
\vspace{0.5ex}
$
{\frak G}_\tau^{}(x)=\sum\limits_{i=1}^{n}{\rm g}_{\tau i}^{}(x)\;\!\partial_{x_i^{}}^{}
$
for all 
$
x\in\Omega,$
$\tau =q+1,\ldots, n$
\vspace{0.35ex}
{\rm(}the set of operators ${\frak G}_\tau^{}$ and {\rm (14)} is commutative  
and these operators aren't holomorphic linearly bound 
\vspace{0.25ex}
on the domain $\Omega\subset \K^n),$
and the holomorphic on the domain $\Omega^{n-q}$ from the space $\K^{n-q}$ functions 
$\varphi_{j\;\!i}^{},\ j = 1,\ldots, m,\ i =1,\ldots,n,$ such that
\vspace{0.25ex}
the conditions of completely solvability on the domain $\Omega$ for the system {\rm(16)} hold.
} 
\vspace{1ex}

The set of the linear differential ope\-ra\-tors of first order
\\[1.75ex]
\mbox{}\hfill
$
{\frak G}_{\tau}^{}(x)=x_{\tau}\;\!\partial_{x_{\tau}^{}}
$
\ for all 
$
x\in\K^n,
\quad
\tau=q+1,\ldots,n,
\hfill
$
\\[1.75ex]
and the infinitesimal operators (14) is commutative  
and these operators aren't holomorphic linearly bound 
on the space $\K^n.$
Then, from Proposition 1, we obtain  
\vspace{1ex}

{\bf Proposition 2.}
{\it 
An system of equations in total differentials 
\\[2ex]
\mbox{}\hfill       
$
\displaystyle
dx_i^{} = 
\sum_{j=1}^{m}\;\!
x_i^{} \varphi_{ji}^{}(x_{q+1}^{},\ldots,x_{n}^{})\;\! dt_j^{},
\ \ i = 1,\ldots, n,
\hfill
$
\\[2ex]
where
\vspace{0.35ex}
$\varphi_{ji}^{},\ j=1,\ldots, m, \ i=1,\ldots,n,$
are holomorphic functions on a domain $\Omega^{n-q}$ from the space $\K^{n-q}$
such that the Frobenius conditions hold
\\[2ex]
\mbox{}\hfill
$
\displaystyle
\sum_{\nu=1}^{n}
x_\nu^{}\varphi_{\mu\nu}^{} (x_{q+1}^{},\ldots,x_{n}^{})\;\!
\partial_{x_{\nu}^{}}^{}\bigl(x_i^{}\varphi_{ji}^{} (x_{q+1}^{},\ldots,x_{n}^{})\bigr)
=
\hfill
$
\\[2.25ex]
\mbox{}\hfill
$
\displaystyle
=
\sum_{\nu=1}^{n}
x_{\nu}^{}\varphi_{j\nu}^{} (x_{q+1}^{},\ldots,x_{n}^{})\;\!
\partial_{x_{\nu}^{}}^{} \bigl(x_i^{}\varphi_{\mu i}^{} (x_{q+1}^{},\ldots,x_{n}^{})\bigr)
$
\ for all 
$
x\in\Omega,
\quad 
\Omega\subset\K^n,
\hfill
$
\\[2.25ex]
\mbox{}\hfill
$
i=1,\ldots,n,\ \ 
j=1,\ldots,m,\ \
\mu=1,\ldots,m,
\hfill
$
\\[2ex]
admits the $q\!\!$-parameter abelian Lie group of transformations {\rm (13)} with the
infinitesimal operators {\rm (14)} and
the universal invariant {\rm (15)}.
}
\vspace{1ex}

{\bf Example 2.}
The one-parameter Lie group of dilatations of space
\\[1.75ex]
\mbox{}\hfill                % (17)
$
x_i^{}\to\, e^{\;\!\alpha}x_i^{},
\quad
i=1,\ldots, n,
$
\hfill (17)
\\[1.75ex]
has the infinitesimal operator
\\[2ex]
\mbox{}\hfill               % (18)
$
\displaystyle
{\frak G}_1^{}(x)=
\sum_{i=1}^{n}x_i^{}\;\!\partial _{x_i^{}}^{}
$
\ for all 
$
x\in\K^n
$
\hfill (18)
\\[1.75ex]
and the universal invariant
\\[2.25ex]
\mbox{}\hfill                 %(19)
$
I\colon x\to 
\Bigl(\;\!\dfrac{x_1^{}}{x_n^{}}\,,\ldots,\dfrac{x_{n-1}}{x_{n}}\Bigr)
$
\ for all 
$
x\in\Omega,
\quad 
\Omega\subset\K^n.
$
\hfill (19)
\\[2ex]
\indent
By Theorem 2 (under $q=1,\, m=1),$ we have
\vspace{1ex}

{\bf Proposition 3.}
{\it 
An autonomous ordinary differential system of the $n\!$-th order 
admits the one-parameter Lie group of dilatations {\rm (17)}
with the infinitesimal operator {\rm (18)} and the universal invariant {\rm (19)}
if and only if this system has the form
\\[2.25ex]
\mbox{}\hfill       
$
\displaystyle
\dfrac{dx_i^{}}{dt} =
x_i^{}\;\!\varphi_1^{}\Bigl(\;\!\dfrac{x_1^{}}{x_n^{}}\,,\ldots,\dfrac{x_{n-1}}{x_{n}}\Bigr)+
\sum_{\tau=2}^{n}\;\!
{\rm g}_{\tau i}^{}(x)\;\! \varphi_{\tau}^{}\Bigl(\;\!\dfrac{x_1^{}}{x_n^{}}\,,\ldots,\dfrac{x_{n-1}}{x_{n}}\Bigr),
\quad 
i = 1,\ldots, n,
\hfill
$
\\[2.5ex]
where the holomorphic functions 
${\rm g}_{\tau i}^{}\colon  \Omega\to \K,\ 
\tau= 2,\ldots, n,\  i =1,\ldots, n,$
\vspace{0.5ex}
are the coordinates of the linear differential ope\-ra\-tors of first order
\vspace{0.5ex}
${\frak G}_\tau^{}(x)=\sum\limits_{i=1}^{n}{\rm g}_{\tau i}^{}(x)\;\!\partial_{x_i^{}}^{}
$
for all 
$
x\in\Omega,
$
$
\tau =2,\ldots, n
$ 
\vspace{0.25ex}
{\rm(}the set of operators ${\frak G}_\tau^{}$ and {\rm (18)} is commutative  
and these operators aren't holomorphic linearly bound 
\vspace{0.35ex}
on the domain $\Omega\subset \K^n),$
and the functions $\varphi_{i}^{}\colon W\to\K, \ i=1,\ldots,n,$ 
are holomorphic on the domain $W$ from the space $\K^{n-1}.$
}
\vspace{1ex}

The set of the linear differential ope\-ra\-tors of first order
\\[2ex]
\mbox{}\hfill
$
{\frak G}_{\tau}^{}(x)=x_{\tau}^{}\;\!\partial_{x_{\tau}^{}}
$
\ for all 
$
x\in\K^n,
\quad
\tau=2,\ldots,n,
\hfill
$
\\[2ex]
and the infinitesimal operators (18) is commutative  
and these operators aren't holomorphic linearly bound 
on the space $\K^n.$
Then, from Proposition 3, we have  
\vspace{1ex}

{\bf Proposition 4.}
{\it 
An autonomous ordinary differential system
\\[2.25ex]
\mbox{}\hfill       
$
\displaystyle
\dfrac{dx_i^{}}{dt} =
x_i^{}\varphi_i^{}\Bigl(\dfrac{x_1^{}}{x_n^{}}\,,\ldots,\dfrac{x_{n-1}}{x_{n}}\Bigr),
\quad 
i = 1,\ldots, n,
\hfill
$
\\[2.5ex]
where the functions $\varphi_{i}^{}\colon W\to\K, \ i=1,\ldots,n,$ 
are holomorphic on the domain $W$ from the space $\K^{n-1},$
admits the one-parameter Lie group of dilatations {\rm (17)}
with the infinitesimal operator {\rm (18)} and the universal invariant {\rm (19)}.
}
\vspace{1ex}

{\bf Example 3.}
Consider the one-parameter Lie group of rotation of real plane
\\[1ex]
\mbox{}\hfill                % (20)
$
u= x\cos \alpha-y\sin \alpha, 
\quad
v= x\sin \alpha+y\cos \alpha
$
\hfill (20)
\\[2.25ex]
with the group parameter $\alpha\in\R.$ 
This Lie group has the infinitesimal operator
\\[2ex]
\mbox{}\hfill               % (21)
$
\displaystyle
{\frak G}(x,y)=
{}-y\;\!\partial _x^{}+x\;\!\partial _y^{}
$
\ for all 
$
(x,y)\in\R^2
$
\hfill (21)
\\[2ex]
and the universal invariant
\\[2ex]
\mbox{}\hfill                 %(22)
$
I\colon (x,y)\to x^2+y^2
$
\ for all 
$
(x,y)\in\R^2.
$
\hfill (22)
\\[2.25ex]
\indent
By Theorem 2 (under $q=1,\, m=1,\, n=2),$ we have
\vspace{1ex}

{\bf Proposition 5.}
{\it 
An autonomous ordinary differential system of second order 
admits the one-parameter Lie group of 
rotation of phase plane {\rm (20)}
with the infinitesimal operator {\rm (21)} and the universal invariant {\rm (22)}
if and only if this system has the form
\\[2.25ex]
\mbox{}\hfill       
$
\displaystyle
\dfrac{dx}{dt} =
{}-y\;\!\varphi_1^{}(x^2+y^2)+a_x^{}(x,y)\;\!\varphi_2^{}(x^2+y^2),
\quad
\dfrac{dy}{dt} =
x\;\!\varphi_1^{}(x^2+y^2)+a_y^{}(x,y)\;\!\varphi_2^{}(x^2+y^2),
\hfill
$
\\[2.75ex]
where the holomorphic functions $a_x^{}\colon  \Omega\to \R$ and $a_y^{}\colon  \Omega\to \R$
\vspace{0.5ex}
are the coordinates of the linear differential ope\-ra\-tor of first order
\vspace{0.5ex}
${\frak A}(x,y)=a_x^{}(x,y)\;\!\partial_x^{}+a_y^{}(x,y)\;\!\partial_y^{}$
for all $(x,y)\in\Omega$ 
{\rm(}the operator ${\frak A}$ is an differential operator such that 
\vspace{0.25ex}
${\frak A}$ is commutative with the infinitesimal operator {\rm (21)} 
and these operators aren't holomorphic linearly bound 
\vspace{0.25ex}
on the domain $\Omega\subset \R^2),$
and sections of functions 
\vspace{0.25ex}
$\varphi_{1}^{}\colon [0;{}+\infty)\to\R$ and $\varphi_{2}^{}\colon [0;{}+\infty)\to\R$ 
are holomorphic on a set $W\subset [0;{}+\infty).$
}
\vspace{1ex}

The linear differential ope\-ra\-tor of first order
\\[2ex]
\mbox{}\hfill
$
{\frak A}(x,y)=x\;\!\partial_x^{}+y\;\!\partial_y^{}
$
\ for all 
$
(x,y)\in\R^2
\hfill
$
\\[2ex]
and the infinitesimal operator (21) aren't holomorphic linearly bound on the plane $\R^2$ and 
they are commutative operators. Then, using Proposition 5, we get
\vspace{0.75ex}

{\bf Proposition 6.}
{\it 
An autonomous ordinary differential system
\\[2.25ex]
\mbox{}\hfill                        % (23)
$
\displaystyle
\dfrac{dx}{dt} =
{}-y\;\!\varphi_1^{}(x^2+y^2)+x\;\!\varphi_2^{}(x^2+y^2),
\quad
\dfrac{dy}{dt} =
x\;\!\varphi_1^{}(x^2+y^2)+y\;\!\varphi_2^{}(x^2+y^2),
$
\hfill {\rm (23)}
\\[2.75ex]
where sections of functions 
\vspace{0.35ex}
$\varphi_{1}^{}\colon [0;{}+\infty)\to\R$ and $\varphi_{2}^{}\colon [0;{}+\infty)\to\R$ 
are holomorphic on $W\subset [0;{}+\infty),$
\vspace{0.25ex}
admits the one-parameter Lie group of 
rotation of phase plane {\rm (20)}
with the infinitesimal operator {\rm (21)} and 
the universal invariant {\rm (22)}.
}
\vspace{1ex}

Let us remark that a differential equation of the first order that
admit the Lie group of rotation {\rm (20)} was considered in [13, p. 149].
This equation is the equation of trajectories for the 
autonomous ordinary differential system (23).
\vspace{0.35ex}

There exists (Theorem 4.1 in [14, p. 23]) the formal change of variables
\\[2ex]
\mbox{}\hfill                        % (24)
$
\displaystyle
u=x+\sum_{k=2}^{\infty}\, U_k^{}(x,y),
\qquad 
v=y+\sum_{k=2}^{\infty}\, V_k^{}(x,y),
$
\hfill {\rm (24)}
\\[2.25ex]
where  
\vspace{0.75ex}
$U_k^{}\colon \R^2\to\R$ and $V_k^{}\colon \R^2\to\R$ are homogeneous polynomials of degrees 
${\rm deg}\, U_k^{}(x,y)=
\linebreak
={\rm deg}\, V_k^{}(x,y)=k,
\ k=2, 3,\ldots,$  
such that the differential system
\\[2.25ex]
\mbox{}\hfill                        % (25)
$
\displaystyle
\dfrac{du}{dt} ={}-v-\sum_{k=2}^{\infty}\, P_k^{}(u,v),
\qquad 
\dfrac{dv}{dt} =u+\sum_{k=2}^{\infty}\, Q_k^{}(u,v),
$
\hfill {\rm (25)}
\\[2.5ex]
where 
\vspace{0.75ex}
$P_k^{}\colon \R^2\!\to\R\!$ and $\!Q_k^{}\colon \R^2\!\to\R$ 
are homogeneous polynomials of degrees ${\rm deg}\, P_k^{}(u,v)=
\linebreak
={\rm deg}\, Q_k^{}(u,v)=k,\ k=2, 3,\ldots,$ 
reduces to the system (23) with $\varphi_1^{}(0)=1,\ \varphi_2^{}(0)=0.$  
\vspace{0.75ex}

Therefore, we have 
\vspace{0.5ex}

{\bf Proposition 7.}
{\it 
Suppose the autonomous ordinary differential system {\rm(25)} has 
an equilibrum state with purely imaginary characteristic roots. Then, 
there exists the formal change of dependent variables {\rm (24)}
such that the system {\rm(25)} reduces to an autonomous differential system that 
admit the one-parameter Lie group of rotation of phase plane.
}
\vspace{0.75ex}

The autonomous differential system (23) that admit 
the one-parameter Lie group of rotation of phase plane (20) is 
the normal form [15] of the autonomous differential system (25) with
an equilibrum state as purely imaginary characteristic roots. 
\vspace{0.75ex}

{\bf Example 4.}
Consider the one-parameter Lie group of Lorentz transformations of real plane
\\[1.75ex]
\mbox{}\hfill                % (26)
$
u= x\cosh \alpha+y\sinh \alpha, 
\quad
v= x\sinh \alpha+y\cosh \alpha
$
\hfill (26)
\\[2.25ex]
with the group parameter $\alpha\in\R.$ 
This Lie group has the infinitesimal operator
\\[2ex]
\mbox{}\hfill               % (27)
$
\displaystyle
{\frak G}(x,y)=
y\;\!\partial _x^{}+x\;\!\partial _y^{}
$
\ for all 
$
(x,y)\in\R^2
$
\hfill (27)
\\[2ex]
and the universal invariant
\\[2ex]
\mbox{}\hfill                 %(28)
$
I\colon (x,y)\to y^2-x^2
$
\ for all 
$
(x,y)\in\R^2.
$
\hfill (28)
\\[2.25ex]
\indent
By Theorem 2 (under $q=1,\, m=1,\, n=2),$ we have
\vspace{1ex}

{\bf Proposition 8.}
{\it 
An autonomous ordinary differential system of second order 
admits the one-parameter Lie group of 
Lorentz transformations of phase plane {\rm (26)}
with the infinitesimal operator {\rm (27)} and the universal invariant {\rm (28)}
if and only if this system has the form
\\[2.25ex]
\mbox{}\hfill       
$
\displaystyle
\dfrac{dx}{dt} =
y\;\!\varphi_1^{}(y^2-x^2)+a_x^{}(x,y)\;\!\varphi_2^{}(y^2-x^2),
\quad
\dfrac{dy}{dt} =
x\;\!\varphi_1^{}(y^2-x^2)+a_y^{}(x,y)\;\!\varphi_2^{}(y^2-x^2),
\hfill
$
\\[2.75ex]
where the holomorphic functions $a_x^{}\colon  \Omega\to \R$ and $a_y^{}\colon  \Omega\to \R$
\vspace{0.5ex}
are the coordinates of the linear differential ope\-ra\-tor of first order
\vspace{0.5ex}
${\frak A}(x,y)=a_x^{}(x,y)\;\!\partial_x^{}+a_y^{}(x,y)\;\!\partial_y^{}$
for all $(x,y)\in\Omega$ 
{\rm(}the operator ${\frak A}$ is an differential operator such that 
\vspace{0.25ex}
${\frak A}$ is commutative with the infinitesimal operator {\rm (27)} 
and these operators aren't holomorphic linearly bound 
\vspace{0.25ex}
on the domain $\Omega\subset \R^2),$
and functions $\varphi_{1}^{}\colon W\to\R$ and $\varphi_{2}^{}\colon W\to\R$ 
are holomorphic on a set $W\subset \R.$
}
\vspace{1ex}

The linear differential ope\-ra\-tor of first order
\\[2ex]
\mbox{}\hfill
$
{\frak A}(x,y)=x\;\!\partial_x^{}+y\;\!\partial_y^{}
$
\ for all 
$
(x,y)\in\R^2
\hfill
$
\\[2ex]
and the infinitesimal operator (27) aren't holomorphic linearly bound on the plane $\R^2$ and 
they are commutative operators. Then, using Proposition 8, we get
\vspace{0.75ex}

{\bf Proposition 9.}
{\it 
An autonomous ordinary differential system
\\[2.25ex]
\mbox{}\hfill       
$
\displaystyle
\dfrac{dx}{dt} =
y\;\!\varphi_1^{}(y^2-x^2)+x\;\!\varphi_2^{}(y^2-x^2),
\quad
\dfrac{dy}{dt} =
x\;\!\varphi_1^{}(y^2-x^2)+y\;\!\varphi_2^{}(y^2-x^2),
\hfill
$
\\[2.75ex]
where functions $\varphi_{1}^{}\colon W\to\R$ and $\varphi_{2}^{}\colon W\to\R$ 
\vspace{0.25ex}
are holomorphic on $W\subset \R,$
admits the one-parameter Lie group of 
Lorentz transformations of phase plane {\rm (26)}
with the infinitesimal operator {\rm (27)} and 
the universal invariant {\rm (28)}.
}
\vspace{1ex}

{\bf Example 5.}
Consider the one-parameter Lie group of projective transformations of real plane
\\[1ex]
\mbox{}\hfill                % (29)
$
u= \dfrac{x}{1-\alpha x}\,, 
\qquad
v= \dfrac{y}{1-\alpha x}
$
\hfill (29)
\\[2.25ex]
with the group parameter $\alpha\in\R.$ 
This Lie group has the infinitesimal operator
\\[2ex]
\mbox{}\hfill               % (30)
$
\displaystyle
{\frak G}(x,y)=
x^2\;\!\partial _x^{}+xy\;\!\partial _y^{}
$
\ for all 
$
(x,y)\in\R^2
$
\hfill (30)
\\[2ex]
and the universal invariant
\\[2.25ex]
\mbox{}\hfill                 %(31)
$
I\colon (x,y)\to \dfrac{x}{y}
$
\, \ for all 
$
(x,y)\in \{ (x,y)\colon y\ne 0\}.
$
\hfill (31)
\\[2.15ex]
\indent
By Theorem 2 (under $q=1,\, m=1,\, n=2),$ we get
\vspace{1ex}

{\bf Proposition 10.}
{\it 
An autonomous ordinary differential system of second order 
admits the one-parameter Lie group of 
projective transformations of phase plane {\rm (29)}
with the infinitesimal operator {\rm (30)} and the universal invariant {\rm (31)}
if and only if this system has the form
\\[2.25ex]
\mbox{}\hfill       
$
\displaystyle
\dfrac{dx}{dt} =
x^2\,\varphi_1^{}\Bigl(\dfrac{x}{y}\Bigr)+a_x^{}(x,y)\;\!\varphi_2^{}\Bigl(\dfrac{x}{y}\Bigr),
\qquad
\dfrac{dy}{dt} =
xy\;\!\varphi_1^{}\Bigl(\dfrac{x}{y}\Bigr)+a_y^{}(x,y)\;\!\varphi_2^{}\Bigl(\dfrac{x}{y}\Bigr),
\hfill
$
\\[2.75ex]
where the holomorphic functions $a_x^{}\colon  \Omega\to \R$ and $a_y^{}\colon  \Omega\to \R$
\vspace{0.5ex}
are the coordinates of the linear differential ope\-ra\-tor of first order
\vspace{0.5ex}
${\frak A}(x,y)=a_x^{}(x,y)\;\!\partial_x^{}+a_y^{}(x,y)\;\!\partial_y^{}$
for all $(x,y)\in\Omega$ 
{\rm(}the operator ${\frak A}$ is an differential operator such that 
\vspace{0.25ex}
${\frak A}$ is commutative with the infinitesimal operator {\rm (30)} 
and these operators aren't holomorphic linearly bound 
\vspace{0.25ex}
on the domain $\Omega\subset \R^2),$
and functions $\varphi_{1}^{}\colon W\to\R$ and $\varphi_{2}^{}\colon W\to\R$ 
are holomorphic on a set $W\subset \R.$
}
\vspace{1ex}

The linear differential ope\-ra\-tor of first order
\\[2ex]
\mbox{}\hfill
$
{\frak A}(x,y)=xy\;\!\partial_x^{}+(x+y^2)\;\!\partial_y^{}
$
\ for all 
$
(x,y)\in\R^2
\hfill
$
\\[2ex]
and the infinitesimal operator (30) aren't holomorphic linearly bound on the plane $\R^2$ and 
they are commutative operators. Then, using Proposition 10, we have
\vspace{0.75ex}

{\bf Proposition 11.}
{\it 
An autonomous ordinary differential system
\\[2.25ex]
\mbox{}\hfill       
$
\displaystyle
\dfrac{dx}{dt} =
x^2\,\varphi_1^{}\Bigl(\dfrac{x}{y}\Bigr)+xy\;\!\varphi_2^{}\Bigl(\dfrac{x}{y}\Bigr),
\qquad
\dfrac{dy}{dt} =
xy\;\!\varphi_1^{}\Bigl(\dfrac{x}{y}\Bigr)+(x+y^2)\;\!\varphi_2^{}\Bigl(\dfrac{x}{y}\Bigr),
\hfill
$
\\[2.75ex]
where functions $\varphi_{1}^{}\colon W\to\R$ and $\varphi_{2}^{}\colon W\to\R$ 
\vspace{0.25ex}
are holomorphic on $W\subset \R,$
admits the one-parameter Lie group of 
projective transformations of phase plane {\rm (29)}
with the infinitesimal operator {\rm (30)} and 
the universal invariant {\rm (31)}.
}
\vspace{1ex}

{\bf Example 6.}
Consider the one-parameter Lie group of nonhomogeneous stretches of real plane
\\[0.75ex]
\mbox{}\hfill                % (32)
$
u= e^{\;\!\alpha}\;\! x, 
\qquad
v= e^{\;\! k \alpha}\;\! y
$
\hfill (32)
\\[2ex]
with the group parameter $\alpha\in\R$ and the coefficient $k\in \R.$ 
This Lie group has the infinitesimal operator
\\[1.5ex]
\mbox{}\hfill               % (33)
$
\displaystyle
{\frak G}(x,y)=
x\;\!\partial _x^{}+ky\;\!\partial _y^{}
$
\ for all 
$
(x,y)\in\R^2
$
\hfill (33)
\\[2ex]
and the universal invariant
\\[2ex]
\mbox{}\hfill                 %(34)
$
I\colon (x,y)\to\ \dfrac{x^k}{y}
$
\ \, for all 
$
(x,y)\in \{ (x,y)\colon y\ne 0\}.
$
\hfill (34)
\\[2.15ex]
\indent
By Theorem 2 (under $q=1,\, m=1,\, n=2),$ we obtain
\vspace{1.25ex}

{\bf Proposition 12.}
{\it 
An autonomous ordinary differential system of second order 
admits the one-parameter Lie group of 
nonhomogeneous stretches of phase plane {\rm (32)}
with the infinitesimal operator {\rm (33)} and the universal invariant {\rm (34)}
if and only if this system has the form
\\[2.25ex]
\mbox{}\hfill       
$
\displaystyle
\dfrac{dx}{dt} =
x\;\!\varphi_1^{}\Bigl(\dfrac{x^k}{y}\Bigr)+a_x^{}(x,y)\;\!\varphi_2^{}\Bigl(\dfrac{x^k}{y}\Bigr),
\qquad
\dfrac{dy}{dt} =
ky\;\!\varphi_1^{}\Bigl(\dfrac{x^k}{y}\Bigr)+a_y^{}(x,y)\;\!\varphi_2^{}\Bigl(\dfrac{x^k}{y}\Bigr),
\hfill
$
\\[2.75ex]
where the holomorphic functions $a_x^{}\colon  \Omega\to \R$ and $a_y^{}\colon  \Omega\to \R$
\vspace{0.5ex}
are the coordinates of the linear differential ope\-ra\-tor of first order
\vspace{0.5ex}
${\frak A}(x,y)=a_x^{}(x,y)\;\!\partial_x^{}+a_y^{}(x,y)\;\!\partial_y^{}$
for all $(x,y)\in\Omega$ 
{\rm(}the operator ${\frak A}$ is an differential operator such that 
\vspace{0.25ex}
${\frak A}$ is commutative with the infinitesimal operator {\rm (33)} 
and these operators aren't holomorphic linearly bound 
\vspace{0.25ex}
on the domain $\Omega\subset \R^2),$
and functions $\varphi_{1}^{}\colon W\to\R$ and $\varphi_{2}^{}\colon W\to\R$ 
are holomorphic on a set $W\subset \R.$
}
\vspace{1ex}

The linear differential ope\-ra\-tor of first order
\\[2ex]
\mbox{}\hfill
$
{\frak A}(x,y)=kx\;\!\partial_x^{}+y\;\!\partial_y^{}
$
\ for all 
$
(x,y)\in\R^2
\hfill
$
\\[2ex]
and the infinitesimal operator (33) aren't holomorphic linearly bound on the plane $\R^2$ and 
they are commutative operators. Then, using Proposition 12, we obtain
\vspace{0.75ex}

{\bf Proposition 13.}
{\it 
An autonomous ordinary differential system
\\[1.75ex]
\mbox{}\hfill       
$
\displaystyle
\dfrac{dx}{dt} =
x\;\!\varphi_1^{}\Bigl(\dfrac{x^k}{y}\Bigr)+kx\;\!\varphi_2^{}\Bigl(\dfrac{x^k}{y}\Bigr),
\qquad
\dfrac{dy}{dt} =
ky\;\!\varphi_1^{}\Bigl(\dfrac{x^k}{y}\Bigr)+y\;\!\varphi_2^{}\Bigl(\dfrac{x^k}{y}\Bigr),
\hfill
$
\\[2.15ex]
where functions $\varphi_{1}^{}\colon W\to\R$ and $\varphi_{2}^{}\colon W\to\R$ 
are holomorphic on $W\subset \R,$
admits the one-parameter Lie group of 
nonhomogeneous stretches of phase plane {\rm (32)}
with the infinitesimal operator {\rm (33)} and 
the universal invariant {\rm (34)}.
}
\vspace{1ex}

{\bf Example\! 7.}\!
Consider\! the\! one-parameter\! Lie\! group\! of\! Galilean\! transformations\! of\! real\! plane
\\[1.5ex]
\mbox{}\hfill                % (35)
$
u=x+\alpha\;\! y, 
\ \
v=  y
$
\hfill (35)
\\[1.5ex]
with the group parameter $\alpha\in\R.$ 
This Lie group has the infinitesimal operator
\\[2ex]
\mbox{}\hfill               % (36)
$
\displaystyle
{\frak G}(x,y)=
y\;\!\partial _x^{}
$
\ for all 
$
(x,y)\in\R^2
$
\hfill (36)
\\[1.75ex]
and the universal invariant
\\[1.75ex]
\mbox{}\hfill                 %(37)
$
I\colon (x,y)\to\ y
$
\ for all 
$
(x,y)\in \R^2.
$
\hfill (37)
\\[1.75ex]
\indent
By Theorem 2 (under $q=1,\, m=1,\, n=2),$ we have
\vspace{0.75ex}

{\bf Proposition 14.}
{\it 
An autonomous ordinary differential system of second order 
admits the one-parameter Lie group of 
Galilean transformations of phase plane {\rm (35)}
with the infinitesimal operator {\rm (36)} and the universal invariant {\rm (37)}
if and only if this system has the form
\\[2ex]
\mbox{}\hfill       
$
\displaystyle
\dfrac{dx}{dt} =
y\;\!\varphi_1^{}(y)+a_x^{}(x,y)\;\!\varphi_2^{}(y),
\qquad
\dfrac{dy}{dt} =a_y^{}(x,y)\;\!\varphi_2^{}(y),
\hfill
$
\\[2.25ex]
where the holomorphic functions $a_x^{}\colon  \Omega\to \R$ and $a_y^{}\colon  \Omega\to \R$
\vspace{0.5ex}
are the coordinates of the linear differential ope\-ra\-tor of first order
\vspace{0.5ex}
${\frak A}(x,y)=a_x^{}(x,y)\;\!\partial_x^{}+a_y^{}(x,y)\;\!\partial_y^{}$
for all $(x,y)\in\Omega$ 
{\rm(}the operator ${\frak A}$ is an differential operator such that 
\vspace{0.25ex}
${\frak A}$ is commutative with the infinitesimal operator {\rm (36)} 
and these operators aren't holomorphic linearly bound 
\vspace{0.25ex}
on the domain $\Omega\subset \R^2),$
and functions $\varphi_{1}^{}\colon W\to\R$ and $\varphi_{2}^{}\colon W\to\R$ 
are holomorphic on a set $W\subset \R.$
}
\vspace{1ex}

The linear differential ope\-ra\-tor of first order
\\[1.75ex]
\mbox{}\hfill
$
{\frak A}(x,y)=x\;\!\partial_x^{}+y\;\!\partial_y^{}
$
\ for all 
$
(x,y)\in\R^2
\hfill
$
\\[2ex]
and the infinitesimal operator (36) aren't holomorphic linearly bound on the plane $\R^2$ and 
they are commutative operators. Then, using Proposition 14, we obtain
\vspace{0.75ex}

{\bf Proposition 15.}
{\it 
An autonomous ordinary differential system
\\[1.75ex]
\mbox{}\hfill       
$
\displaystyle
\dfrac{dx}{dt} =
y\;\!\varphi_1^{}(y)+x\;\!\varphi_2^{}(y),
\qquad
\dfrac{dy}{dt} =y\;\!\varphi_2^{}(y),
\hfill
$
\\[2ex]
where functions $\varphi_{1}^{}\colon W\to\R$ and $\varphi_{2}^{}\colon W\to\R$ 
are holomorphic on $W\subset \R,$
admits the one-parameter Lie group of 
Galilean transformations of phase plane {\rm (35)}
with the infinitesimal operator {\rm (36)} and 
the universal invariant {\rm (37)}.
}
\vspace{1.5ex}

\newpage

\mbox{}
\\[-3.5ex]

}
\end{document}